\documentclass[11pt]{article}

\usepackage{fullpage}
\usepackage{amsgen}
\usepackage{amsmath}
\usepackage{amstext}
\usepackage{amsbsy}
\usepackage{amsopn}
\usepackage{amsfonts}
\usepackage{amssymb}
\usepackage{eepic}
\usepackage{graphicx}
\usepackage{epsf}
\usepackage{pstricks}
\usepackage{pstricks-add}
\usepackage{pgf,tikz}

\def\Box{\square}

\def\mapright#1{\smash{\mathop{\longrightarrow}\limits^{#1}}}
\def\tra#1{\smash{\mathop{\mid\kern
-1pt\joinrel\relbar\joinrel\relbar}\limits^{*}_{#1}}}
\def\longtra#1{\smash{\mathop{\mid\kern
-1pt\joinrel\relbar\joinrel\relbar\joinrel\relbar}\limits^{*}_{#1}}}
\def\vlongtra#1{\smash{\mathop{\mid\kern
-1pt\joinrel\relbar\joinrel\relbar\joinrel\relbar\joinrel\relbar}\limits^{*}_{#1}}}
\def\vvlongtra#1{\smash{\mathop{\mid\kern
-1pt\joinrel\relbar\joinrel\relbar\joinrel\relbar\joinrel\relbar\joinrel\relbar}\limits^{*}_{#1}}}
\def\vvvlongtra#1{\smash{\mathop{\mid\kern
-1pt\joinrel\relbar\joinrel\relbar\joinrel\relbar\joinrel\relbar\joinrel\relbar\joinrel\relbar}\limits^{*}_{#1}}}
\def\etra#1{\smash{\mathop{\mid\kern
-1pt\joinrel\relbar\joinrel\relbar}\limits_{#1}}}

\def\iff{\Leftrightarrow}
\def\Rw{\Rightarrow}
\def\oo{\overline}

\def\wt{\widetilde}

\def\mesh{\mbox{mesh}}

\def\im{\mbox{Im}}
\def\haus{\mbox{Haus}}

\def\N{\mathbb{N}}
\def\RR{\mathbb{R}}

\def\diam{\mbox{diam}}

\def\dom{\mbox{dom}}

\def\max{\mbox{max}}
\def\min{\mbox{min}}

\def\sup{\mbox{sup}}

\def\p{\varphi}

\def\inv{^{-1}}

\def\bi{\begin{itemize}}
\def\ei{\end{itemize}}
\def\beq{\begin{equation}}
\def\eeq{\end{equation}}

\newtheorem{T}{Theorem}[section]
\newcommand{\bt}{\begin{T}}
\newcommand{\et}{\end{T}}
\newcommand{\ftd}{$\square$\end{T}}

\newtheorem{Proposition}[T]{Proposition}
\newcommand{\bp}{\begin{Proposition}}
\newcommand{\ep}{\end{Proposition}}
\newcommand{\fpd}{$\square$\end{Proposition}}

\newtheorem{Lemma}[T]{Lemma}
\newcommand{\bl}{\begin{Lemma}}
\newcommand{\el}{\end{Lemma}}
\newcommand{\fld}{$\square$\end{Lemma}}

\newtheorem{Corol}[T]{Corollary}
\newcommand{\bc}{\begin{Corol}}
\newcommand{\ec}{\end{Corol}}
\newcommand{\fcd}{$\square$\end{Corol}}

\newtheorem{Remark}[T]{Remark}
\newcommand{\br}{\begin{Remark}}
\newcommand{\er}{\end{Remark}}
\newcommand{\frd}{$\square$\end{Remark}}

\newtheorem{Example}[T]{Example}
\newcommand{\be}{\begin{Example}}
\newcommand{\ee}{\end{Example}}

\newtheorem{Problem}[T]{Problem}
\newcommand{\bq}{\begin{Problem}}
\newcommand{\eq}{\end{Problem}}

\newcommand{\proof}
   {\par\medbreak\noindent{\bf Proof}.\enspace}

\newcommand{\qed}{%\hfill
$\Box$
\par\bigbreak}

\normalbaselines
\def\abstract#1{\par\bigskip
\begingroup\small
\baselineskip=12truept
\begin{center}ABSTRACT\end{center}
\par\medskip\par\noindent
\null\hfill\hbox{\vbox{\hsize=5truein\noindent#1}}
\hfill\null\par\endgroup\par}

\title{Geometric characterizations of virtually free groups}
\author{{\bf V\'\i tor Ara\'ujo}\\
$ $\\ {\em Universidade Federal da Bahia, Instituto de
  Matem\'atica,}\\
{\em Av. Adhemar de Barros, S/N, Ondina,}\\
{\em 40170-110 Salvador-BA, Brazil}\\
{\em email:} vitor.d.araujo@ufba.br\\
$ $\\
{\bf Pedro V. Silva}\\ $ $\\
{\em Universidade Federal da Bahia, Instituto de
  Matem\'atica, Brazil}\\
{\em and Centro de
Matem\'{a}tica, Faculdade de Ci\^{e}ncias, Universidade do
Porto,}\\ {\em R. Campo Alegre 687, 4169-007 Porto, Portugal}\\
{\em email:} pvsilva@fc.up.pt}

\date{\today}

\begin{document}
\maketitle

\begin{center}\small
2010 Mathematics Subject Classification: 20F67, 51M10, 54E35

\bigskip

Keywords: hyperbolic metric spaces, hyperbolic groups, virtually free groups 
\end{center}

\abstract{Four geometric conditions on a geodesic metric space, which
  are stronger variants of classical 
  conditions characterizing hyperbolicity, are
  proved to be equivalent. In the particular case of the Cayley graph
  of a finitely generated group, it is shown that they characterize
  virtually free groups.}

\section{Introduction}

Finitely generated virtually free groups constitute an important
subclass of hyperbolic groups. A group $G$ is {\em virtually free} if it
has a free subgroup $F$ of finite index. If $G$
is finitely generated, we may assume the same for $F$.

Virtually free groups constitute probably the class of groups which admits the
widest variety of characterizations. Following Diekert and Weiss
\cite{DW,DW2}, we enumerate a few to indulge all tastes.

A finitely
generated group $G = \langle A\rangle$ is virtually free if and only if:
\bi
\item
$G$ is the fundamental group of a finite graph of finite groups
\cite{KPS2} (see also \cite[Theorem 7.3]{SW});
\item
$G$ acts on a connected locally finite graph of finite treewidth, with
finitely many orbits and finite node stabilizers \cite{DW};
\item
the Cayley graph $\Gamma_A(G)$ has finite treewidth \cite{KL};
\item
$\Gamma_A(G)$ can be $k$-triangulated \cite{MS2};
\item
there exists some $\varepsilon \geq 0$ such that, for all
coterminal geodesic $\xi$ 
and path $\xi'$ in $\Gamma_A(G)$, $\im(\xi) \subseteq
D_{\varepsilon}(\im(\xi'))$ \cite{Ant};
\item
there exist some finite
generating set $B$ of $G$ and some $k \geq 0$
such that every $k$-locally geodesic in $\Gamma_B(G)$ is a geodesic
\cite{GHHR};
\item
$G$ admits a finite presentation by some geodesic rewriting system
\cite{GHHR};
\item
the language of all words on $A \cup A\inv$ representing the identity is 
context-free (Muller-Schupp's Theorem \cite{MS,Dun});
\item
$G$ is the universal group of a finite pregroup \cite{Rim};
\item
the monadic second-order theory of $\Gamma_A(G)$ is decidable \cite{KL,MS2};
\item
$G$ admits a Stallings section \cite{SSV}.
\ei
More details can be found in \cite{Ant,DW2}. 

The above list gives evidence of the interest which was devoted by
different authors to the Cayley graph of a virtually free
group. However, to our best knowledge, no explicit results were
published on the classical geometric conditions used 
to define hyperbolicity \cite{Gro}, \cite[Proposition
2.21]{GH}. 

Thus in Section 3 we establish equivalent geometric
conditions for geodesic metric spaces which will
allow us to characterize
finitely generated virtually free groups when we consider their Cayley
graphs in Section 4. 

Our starting point were three of the most common characterizations of
hyperbolicity, using respectively thin geodesic triangles (Rips condition),
the Gromov product (the original definition by Gromov) and the mesh of
a geodesic triangle. We succeed on replacing thin geodesic triangles
by thin geodesic polygons, and the mesh of
a geodesic triangle by the mesh of an arbitrary triangle. For the
Gromov product, we consider inequalities with arbitrary long sequences
of group elements.
% Thus our approach in this paper is to establish equivalent geometric
% conditions for geodesic metric spaces which allow us to characterize
% finitely generated virtually free groups when we consider their Cayley
% graphs. These tasks are performed in Sections 3 and 4, respectively.

\section{Preliminaries}
\label{hg}

We present in this section well-known facts regarding hyperbolic
spaces and hyperbolic groups. The reader is referred to \cite{BH,GH} for
details.

If $(X,d)$ is a metric space and $\varepsilon \geq 0$, we denote by
$D_{\varepsilon}(x)$ the closed ball of center $x \in X$ and radius
$\varepsilon$. If $Y \subseteq X$ is nonempty, we write
$$D_{\varepsilon}(Y) = \bigcup_{y \in Y} D_{\varepsilon}(y).$$

A mapping $\p:(X,d) \to (X',d')$ between metric spaces is called
% \bi
% \item
% a {\em contraction} if $d'(x\p,y\p) \leq d(x,y)$ for all $x,y
% \in X$;
% \item
an {\em isometric embedding} if $d'(x\p,y\p) = d(x,y)$ for all $x,y
\in X$. A surjective isometric embedding is an
{\em isometry}.

We consider the usual metric for $\RR$ and its subsets. A {\em path}
in $(X,d)$ is a continuous mapping $\xi:[0,s] \to X$ $(s \in
\RR_0^+)$. If $\xi':[0,s'] \to X$ is another path and $0\xi' = 0\xi$,
$s'\xi' = s\xi$, we say that the paths $\xi$ and $\xi'$ are {\em
  coterminal}.

A metric space $(X,d)$ is said to be {\em geodesic} if, for all $x,y
\in X$, there exists some path $\xi:[0,s] \to X$ which is an isometric
embedding and such that $0\xi = x$ and $s\xi = y$. 
We call $\xi$ a {\em geodesic} of
$(X,d)$. 
% We shall often call $\im(\xi)$ a geodesic as well. In this
% second sense, we may use the notation $[x,y]$ to denote an arbitrary geodesic
% connecting $x$ to $y$.

A {\em quasi-isometric embedding} of metric spaces is a mapping $\p:(X,d) \to
(X',d')$ such that there exist constants $\lambda \geq 1$ and $K \geq 0$
satisfying
%\beq
%\label{qie}
$$\frac{1}{\lambda}d(x,y) -K \leq d'(x\p,y\p) \leq \lambda d(x,y) +K$$
%\eeq\frac{1}{\lambda}d(x,y) -K \leq d'(x\p,y\p)
for all $x,y \in X$. We may call it a
$(\lambda,K)$-quasi-isometric embedding if we want to stress the
constants. 
If in addition
%\beq
%\label{qie2}
$$\forall x' \in X'\; \exists x \in X:\; d'(x',x\p) \leq K,$$
%\eeq
we say that $\p$ is a {\em quasi-isometry}.

Two metric spaces $(X,d)$ and $(X',d')$ are said to be {\em
  quasi-isometric} if there exists a quasi-isometry $\p:(X,d) \to
(X',d')$. Quasi-isometry turns out to be an equivalence relation on
the class of metric spaces.
A path $\xi:[0,s] \to X$ which is a quasi-isometric
embedding is called a {\em quasi-geodesic} of $(X,d)$.

Let $(X,d)$ be a geodesic metric space. A {\em triangle} in $(X,d)$ is
a collection $\Delta = [[\xi_0,\xi_1,\xi_2]]$ of three paths
$\xi_i:[0,s_i] \to X$ such that 
$$s_i\xi_i = 0\xi_{i+1}\hspace{1cm}\mbox{ for }i = 0,1,2.$$
Here and in many other instances of the paper, we consider the
indices modulo $3$ (or $n+1$), so that $\xi_{3} = \xi_0$ in the
formula above.
If all the paths $\xi_i$ are geodesics,
%(respectively $(\lambda,K)$-quasi-geodesics), 
we say that $\Delta$ is a {\em geodesic triangle}.
% (respectively $(\lambda,K)$-{\em quasi-geodesic triangle}). 

More generally, we define a {\em geodesic polygon} to be a collection of
$n+1$ geodesics
$\Pi = [[\xi_0,\xi_1,\ldots,\xi_n]]$ such that 
$$s_i\xi_i = 0\xi_{i+1}\hspace{1cm}\mbox{ for }i = 0,1,n,$$
with $n\geq 1$ arbitrary. We may call $\Pi$ an $(n+1)$-{\em gon}. If
$n = 1$, we have a {\em geodesic bigon}. 

Given $\delta \geq 0$, a geodesic polygon
$[[\xi_0,\xi_1,\ldots,\xi_n]]$ is said to be $\delta$-{\em thin} if
$$\im(\xi_n) \subseteq D_{\delta}(\im(\xi_0) \cup \im(\xi_{n-1})).$$
We say
that $(X,d)$ is $\delta$-{\em hyperbolic} if every geodesic triangle in 
$(X,d)$ is $\delta$-thin. We say that $(X,d)$ is {\em hyperbolic} if
$(X,d)$ is $\delta$-hyperbolic for some $\delta \geq 0$. 

% Note that if $n \geq 2$ and all $(n+1)$-gons are $\delta$-thin, then
% all $n$-gons are $\delta$-thin as well.

There are several equivalent characterizations of hyperbolicity. The
original one, introduced by Gromov in \cite{Gro}, uses
the concept of Gromov product, which 
we now define. We note that it can be defined for every metric space $(X,d)$.

Given $x,y,p \in X$, we 
define 
$$(x|y)_p = \frac{1}{2}(d(p,x) + d(p,y)-d(x,y)).$$
We say that $(x|y)_p$ is the {\em Gromov product} of $x$ and $y$,
taking $p$ as basepoint. 

A third road to hyperbolicity uses the concept of {\em mesh}. Recall
the notion of {\em diameter}. Given a nonempty $Y \subseteq X$, we
write
$$\diam(Y) = \sup \{ d(y,y') \mid y,y' \in Y\}.$$

Given a
triangle $\Delta = [[\xi_0,\xi_1,\xi_2]]$ in a geodesic metric space
$(X,d)$, we define 
$$\mesh(\Delta) = \inf\{ \diam(\{ u_0,u_1,u_2\}) \mid u_i \in
\im(\xi_i)\; (i = 0,1,2)\}.$$

The following equivalences are well known, see \cite[Proposition
2.21]{GH} for a proof. 

\bt
\label{hyp}
The following
conditions are equivalent for a geodesic metric space $(X,d)$:
\bi
\item[(i)] $(X,d)$ is hyperbolic;
\item[(ii)] there exists some $\delta' \geq 0$ such that
%\beq
%\label{hyp1}
$$(x_0|x_2)_p \geq {\rm min}\{ (x_0|x_1)_p, (x_1|x_2)_p \} -\delta'$$
%\eeq
for all $x_0,x_1,x_2,p \in X$;
\item[(iii)] there exists some $\mu \geq 0$ such that
$${\rm mesh}(\Delta) \leq \mu$$
for every geodesic triangle $\Delta$ in $(X,d)$.
\ei
\et

Among the most important properties of hyperbolic spaces, stands the fact that
they are closed under quasi-isometry.

Given $Y,Z \subseteq X$ nonempty, the {\em Hausdorff distance} between
$Y$ and $Z$ is defined by
$$\haus(Y,Z) = \max\{ \sup_{y\in Y} d(y,Z),\; \sup_{z\in Z} d(z,Y) \}.$$
The following result of Gromov \cite{Gro} (see also \cite[Theorem
5.4.21]{GH}) is also important to us:

\bp
\label{conr}
Let $\lambda \geq 1$ and $K, \delta \geq 0$. There exists a constant
$R(\delta,\lambda,K) \geq 0$ such that: if $(X,d)$ is a $\delta$-hyperbolic
space, $\xi$ is a geodesic
and $\xi'$ a coterminal $(\lambda,K)$-quasi-geodesic in $(X,d)$, then
$${\rm Haus}({\rm Im}(\xi), {\rm Im}(\xi')) \leq
R(\delta,\lambda,K).$$
\ep

We proceed now to define hyperbolic groups.

Given a subset $A$ of a group $G$, we denote by $\langle A\rangle$ the
subgroup of $G$ generated by $A$. We assume throughout the paper that
generating sets are finite. 

Given $G = \langle A\rangle$, we write $\wt{A} =
A \cup A\inv$. The {\em Cayley
  graph} $\Gamma_A(G)$ has vertex set $G$ and edges of the form $g
\mapright{a} ga$ for all $g \in G$ and $a \in \wt{A}$.
The {\em geodesic metric} $d_A$ on $G$ is defined by taking $d_A(g,h)$
to be the length of the shortest path connecting $g$ to $h$ in $\Gamma_A(G)$. 

Since $\im(d_A) \subseteq \N$, then $(G,d_A)$ is not a geodesic metric
space. However, we can remedy that by embedding $(G,d_A)$
isometrically into the {\em geometric 
realization} $\oo{\Gamma}_A(G)$ of $\Gamma_A(G)$, when vertices become
points and edges 
become segments of length 1 in some (euclidean) space, intersections
being determined by adjacency only. With the obvious metric,
$\oo{\Gamma}_A(G)$ is a geodesic metric space, and the geometric 
realization is unique up to isometry. We denote also by $d_A$ the
induced metric on $\oo{\Gamma}_A(G)$.

% Given $X \subseteq G$
% nonempty and $g \in G$, we define
% $$d_A(g,X) = \min\{ d_A(g,x) \mid x \in X\}.$$
We say that the group $G$ is {\em hyperbolic} if the geodesic metric
space $(\oo{\Gamma}_A(G),d_A)$ is
hyperbolic. Since $(\oo{\Gamma}_A(G),d_A)$ is quasi-isometric to
$(\oo{\Gamma}_B(G),d_B)$ for every alternative finite generating set
$B$ of $G$, the concept of hyperbolic group
is independent from 
the finite generating set considered. However, the hyperbolicity
constant $\delta$ may vary with the generating set. 

Let $H$ be a subgroup of a hyperbolic group $G = \langle A\rangle$ and
let $q \geq 0$. We
say that $H$ is $q$-{\em quasi-convex} with respect to $A$ if, for
every geodesic $\xi:[0,s] \to \oo{\Gamma}_A(G)$ with $0\xi, s\xi \in
H$, we have
$$\im(\xi) \subseteq D_q(H).$$
We say that $H$ is {\em quasi-convex} if it is $q$-quasi-convex
for some $q \geq 0$.
Like most other properties in the theory of hyperbolic groups,
quasi-convex does not depend on the finite generating set considered. 

A (finitely generated) subgroup of a hyperbolic group needs not be
hyperbolic, but a quasi-convex subgroup of a hyperbolic group is
always hyperbolic. The converse is not true in general. 
The reader is referred to \cite[Section III.$\Gamma$.3]{BH} for
details on quasi-convex subgroups.

\section{Geodesic metric spaces}

We prove in this section the equivalence of four geometric
conditions on a geodesic metric space.

Given $\delta \geq 0$, we say that a geodesic metric space $(X,d)$ is
{\em polygon $\delta$-hyperbolic} if every geodesic polygon in $(X,d)$ is
$\delta$-thin. We say that $(X,d)$ is
{\em polygon hyperbolic} if it is polygon $\delta$-hyperbolic for
some $\delta \geq 0$. We introduce also the notation
$$\mesh(X,d) = \sup \{ \mesh(\Delta) \mid \Delta \mbox{ is a triangle
  in }(X,d) \}.$$

\bl
\label{cuqi}
Let $\p:(X,d) \to (X',d')$ be a quasi-isometric embedding of geodesic
metric spaces. Then:
\bi
\item[(i)] if $(X',d')$ is polygon hyperbolic, so is $(X,d)$;
\item[(ii)] if ${\rm mesh}(X',d')$ is finite, so is ${\rm mesh}(X,d)$.
\ei
\el

\proof
(i) Assume that $\p$ is a $(\lambda,K)$-quasi-isometry and
$(X',d')$ is polygon $\delta'$-hyperbolic. 
% Then $(X',d')$
% is hyperbolic and since the class of hyperbolic spaces is closed under
% quasi-isometry [???] $(X,d)$ is $\delta$-hyperbolic for some $\delta
% \geq 0$. 
Let $[[\xi_0,\ldots,\xi_n]]$ be a geodesic polygon in 
$(X,d)$. 
%For $i = 0,\ldots,n$, let $\xi_i:[0,s_i] 
%\to [x_{i},x_{i+1}]$ be the geodesic in $[[x_0,\ldots,x_n]]$.

Each geodesic $\xi_i$ induces a $(\lambda,K)$-quasi-geodesic $\xi_i\p$
in $X'$.
%$\xi_i\p:[0,s_i] \to X'$. 
For $i = 0,\ldots,n$, let $\xi'_i$ be a geodesic in $X'$ coterminal
with $\xi_i\p$.
Since $(X',d')$ is $\delta'$-hyperbolic,
it follows from Proposition \ref{conr} that 
\beq
\label{fgvf3}
\haus(\im(\xi_i\p),\im(\xi'_i)) \leq R(\delta,\lambda,K)
\eeq
for $i = 0,\ldots,n$. 

Write $R = R(\delta,\lambda,K)$ and let $y \in \im(\xi_n)$. Then $y\p
\in \im(\xi_n\p)$ and by 
(\ref{fgvf3}) we have $d'(y\p,y') \leq R$ for some $y' \in \im(\xi'_n)$.

Now $[[\xi'_0,\ldots,\xi'_n]]$ is a geodesic polygon in 
$X'$. Since $(X',d')$ is polygon
$\delta'$-hyperbolic, there exist $j \in \{ 0,\ldots,n-1\}$ and $z' \in
\im(\xi'_j)$ such that $d'(y',z') \leq \delta'$. By (\ref{fgvf3}), we have
$d'(z',z\p) \leq R$ for some $z \in \im(\xi_j)$, hence
$$\frac{1}{\lambda}d(y,z) -K \leq d'(y\p,z\p)$$
yields
$$d(y,z) \leq {\lambda}(d'(y\p,z\p)+K) \leq
\lambda(d'(y\p,y')+d'(y',z')+d'(z',z\p)+K) \leq 
\lambda(2R+\delta' +K).$$
Since $y \in \im(\xi_n)$ is arbitrary and $z \in \im(\xi_j)$,
it follows that $(X,d)$ is polygon
$\delta$-hyperbolic for $\delta = \lambda(2R+\delta' +K)$.

(ii) Assume that $\p$ is a $(\lambda,K)$-quasi-isometry and
$\mesh(X',d') = \mu' < +\infty$. Let $[[\xi_0,\xi_1,\xi_2]]$
be a triangle in $(X,d)$. 
% For $i = 0,\ldots,2$, let $\xi_i:[0,s_i] 
% \to [x_{i},x_{i+1}]$ be the geodesic in $[[x_0,\ldots,x_n]]$.
Then $[[\xi_0\p,\xi_1\p,\xi_2\p]]$ is a triangle in $(X',d')$.
Since $\mesh(X',d') = \mu'$, there exist some $u_i \in
\dom(\xi_i\p) = \dom(\xi_i)$, for $i = 0,1,2$, such that
$$\diam(\{ u_0\xi_0\p, u_1\xi_1\p, u_2\xi_2\p \}) \leq \mu'.$$
Since 
$$\frac{1}{\lambda}d(u_i\xi_i,u_j\xi_j) -K \leq
d'(u_i\xi_i\p,u_j\xi_j\p)$$
for all $i,j \in \{ 0,1,2\}$, we get
$$\diam(\{ u_0\xi_0, u_1\xi_1, u_2\xi_2 \}) \leq \lambda(\mu'+K).$$
Therefore 
$$\mesh(X,d) \leq \lambda(\mu'+K) < +\infty.$$
\qed

\bt
\label{eqg}
Let< $(X,d)$ be a geodesic metric space. Then the following
conditions are equivalent:
\bi
\item[(i)] $(X,d)$ is polygon hyperbolic;
\item[(ii)] there exists some $\delta' \geq 0$ such that
$$(x_0|x_n)_p \geq {\rm min}\{ (x_0|x_1)_p, \ldots (x_{n-1}|x_n)_p
\} -\delta'$$
holds for all $n \geq 2$ and  $x_0,\ldots, x_n,p \in X$;
\item[(iii)]
there exists some $\varepsilon \geq 0$ such that, for all coterminal
geodesic $\xi$ 
and path $\xi'$in $(X,d)$, 
$${\rm Im}(\xi) \subseteq D_{\varepsilon}({\rm
  Im}(\xi'));$$ 
\item[(iv)] ${\rm mesh}(X,d) < +\infty$.
\ei
\et

\proof
(i) $\Rw$ (ii). Assume that $(X,d)$ is polygon $\delta$-hyperbolic.
Let $x_0,\ldots, x_n,p
\in X$ with $n \geq 2$. For $i = 0,\ldots,n$, consider a geodesic
$\xi_i:[0,s_i] \to X$ such that $0\xi_i = x_i$ and
$s_i\xi_i = x_{i+1}$. Let $\Pi = [[\xi_0,\ldots,\xi_n]]$ and $Z = 
\im(\xi_0) \cup \ldots \cup \im(\xi_{n-1})$.

Since $(X,d)$ is polygon $\delta$-hyperbolic, it is in particular
$\delta$-hyperbolic.  
%Let $\delta' = 3\delta$. 
By \cite[Lemmas 2.9, 2.31 and 2.32]{Vai}, we have
\beq
\label{fgvf4}
(x_i|x_{i+1})_p \leq d(p,\im(\xi_i)) \leq (x_i|x_{i+1})_p +2\delta
\eeq
for $i = 0,\ldots,n$. Hence
\beq
\label{fgvf5}
\min\{ (x_0|x_1)_p, \ldots (x_{n-1}|x_n)_p \}
\leq \min\{ d(p,\im(\xi_0)), \ldots, d(p,\im(\xi_{n-1})) \}
= d(p,Z).
\eeq
Since $\im(\xi_n)$ is compact,
we have $d(p,\im(\xi_n)) = d(p,q)$ for some $q \in
\im(\xi_n)$. Applying our hypothesis to $\Pi$ we get $d(q,Z) \leq
\delta$ and so 
\beq
\label{fgvf6}
d(p,Z) \leq d(p,q) + d(q,Z) \leq d(p,\im(\xi_n)) +\delta
\leq (x_n|x_{0})_p +3\delta
\eeq
in view of (\ref{fgvf4}).
Now (\ref{fgvf5}) and (\ref{fgvf6}) together yield
$$(x_0|x_{n})_p \geq d(p,Z) -3\delta \geq \min\{ (x_0|x_1)_p,
\ldots (x_{n-1}|x_n)_p \} -3\delta$$
and we are done.

(ii) $\Rw$ (i). Let 
$[[\xi_0,\ldots,\xi_n]]$ 
be a geodesic polygon in 
$(X,d)$ with geodesics $\xi_i:[0,s_i] \to X$ for $i =
0,\ldots,n$. Write $x_i = 0\xi_i = s_{i-1}\xi_{i-1}$ for $i =
0,\ldots,n$.
Let $p \in \im(\xi_n)$ and
let $Z = \im(\xi_0) \cup \ldots \cup \im(\xi_{n-1})$. 

By
\cite[Lemmas 2.9, 2.30 and 2.32]{Vai}, the case $n = 2$ in condition
(ii) implies 
\beq
\label{fgvf7}
(x_i|x_{i+1})_p \leq d(p,\im(\xi_i)) \leq (x_i|x_{i+1})_p
+2\delta'
\eeq
for $i = 0,\ldots,n$.
Hence
$$d(p,Z) = \min\{ d(p,\im(\xi_0)), \ldots, d(p,\im(\xi_{n-1})) \}
\leq \min\{ (x_0|x_{1})_p, \ldots, (x_{n-1}|x_{n})_p \}
+2\delta'.$$
By condition (ii), we get 
$$d(p,Z) \leq (x_0|x_{n})_p +3\delta'.$$
Since $p \in \im(\xi_n)$, it follows from (\ref{fgvf7}) that 
$d(p,Z) \leq 3\delta'$. Therefore $(X,d)$ is polygon
$3\delta'$-hyperbolic.

(i) $\Rw$ (iii). 
Assume that $(X,d)$ is polygon $\delta$-hyperbolic.
Let $\xi:[0,s] \to X$ be a geodesic and $\xi':[0,s'] \to X$ be a
coterminal path. Since $\im(\xi')$ is compact,
we can find finitely many $x_0 = s\xi, x_1, \ldots,x_{n-1}, x_n = 0\xi
\in \im(\xi')$ such 
that $d(x_i,x_{i-1}) \leq 1$ for $i = 1,\ldots,n$.
%$$\im(\xi') \subseteq D_1(x_0) \cup \ldots \cup D_1(x_n).$$
%We may also assume that $x_0 = s\xi$ and $x_n = 0\xi$. 
Now we build a
geodesic polygon $\Pi = [[\xi_0, \ldots, \xi_n]]$, where each $\xi_i$
is a geodesic with endpoints $x_i$ and $x_{i+1}$, and $\xi_n =
\xi$. 

Let $p \in \im(\xi)$.
Since $(X,d)$ is polygon $\delta$-hyperbolic, there exist some
$j \in \{ 0,\ldots,n-1\}$ and $q \in \im(\xi_j)$ such that $d(p,q)
\leq \delta$. Hence 
$$d(p,\im(\xi') \leq d(p,x_j) \leq d(p,q) + d(q,x_j) \leq \delta + 1$$
and so $\im(\xi) \subseteq D_{\delta+1}(\im(\xi'_i))$.

(iii) $\Rw$ (iv).
If we consider a geodesic triangle $[[\xi_0,\xi_1,\xi_2]]$ in $(X,d)$
and view $\im(\xi_0) \cup \im(\xi_1)$ as the image of a single path,
coterminal with $\xi_2$, it follows from condition (iii) that $(X,d)$
is $\varepsilon$-hyperbolic. Therefore, by Theorem \ref{hyp}, 
there exists some $\mu \geq 0$ such that
$$\mesh(\Delta) \leq \mu$$
for every geodesic triangle $\Delta$ in $(X,d)$.

Let $\Delta' = [[\xi'_0,\xi'_1,\xi'_2]]$ be a triangle in
$(X,d)$. Consider a geodesic triangle $\Delta = [[\xi_0,\xi_1,\xi_2]]$
where $\xi_i$ and $\xi'_i$ are coterminal for $i = 0,1,2$. Since
$\mesh(\Delta) \leq \mu$, there exist $p_i \in \im(\xi_i)$ such that 
$$\diam(\{ p_0,p_1,p_2) \}) \leq \mu.$$
Now, also by condition (iii), there exist $q_i \in \im(\xi'_i)$ such
that $d(p_i,q_i) \leq \varepsilon$ for $i = 0,1,2$. Therefore, for all
$i,j \in \{ 0,1,2\}$, we get
$$d(q_i,q_j) \leq d(q_i,p_i) + d(p_i,p_j) + d(p_j,q_j) \leq
2\varepsilon + \mu.$$
Therefore $\mesh(X,d) \leq 2\varepsilon + \mu$ and (iv) holds.

(iv) $\Rw$ (i). 
Assume that $\mesh(X,d) = \mu < +\infty$.
Let 
$\Pi = [[\xi_0,\ldots,\xi_n]]$ 
be a geodesic polygon in 
$(X,d)$ with geodesics $\xi_i:[0,s_i] \to X$ for $i =
0,\ldots,n$. Write $x_i = 0\xi_i = s_{i-1}\xi_{i-1}$ for $i =
0,\ldots,n$.
Let $p \in \im(\xi_n)$, say $p = y\xi_n$, and
let $Z = \im(\xi_0) \cup \ldots \cup \im(\xi_{n-1})$. We prove that
$d(p,Z) \leq \frac{3\mu +1}{2}$.

We may of course assume that $d(p,x_0), d(p,x_n) \geq \frac{\mu+1}{2}$. We can
build three paths $\xi'_0, \xi'_1,\xi'_2$ satisfying the following
conditions:
\bi
\item
$\im(\xi'_0) = [y+\frac{\mu+1}{2},s_n]\xi_n \cup Z$ and has endpoints
$(y+\frac{\mu+1}{2})\xi_n$ and $x_n$;
\item
$\im(\xi'_1) = [0,y-\frac{\mu+1}{2}]\xi_n$ and has endpoints $x_n$ and
$(y-\frac{\mu+1}{2})\xi_n$;
\item
$\im(\xi'_2) = [y-\frac{\mu+1}{2},y+\frac{\mu+1}{2}]\xi_n$ and has
endpoints $(y-\frac{\mu+1}{2})\xi_n$ and 
$(y+\frac{\mu+1}{2})\xi_n$.
\ei
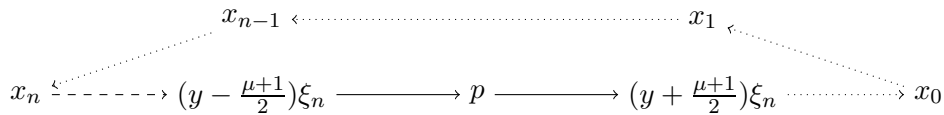
\begin{center}

\medskip

\begin{tikzpicture}[xscale=3,yscale=1]
  \node (a) at (0,0) {$x_n$};
  \node (b) at (1,0) {$(y-\frac{\mu+1}{2})\xi_n$};
  \node (c) at (2,0) {$p$};
  \node (d) at (3,0) {$(y+\frac{\mu+1}{2})\xi_n $};
  \node (e) at (4,0) {$x_0$};
  \node (b1) at (1,1) {$x_{n-1}$};
  \node (d1) at (3,1) {$x_{1}$};
  \draw [->, thin] (b) edge (c);
  \draw [->, thin] (c) edge (d);
  \draw [->, dashed] (a) edge (b);

  \draw [->, dotted] (d1) edge (b1);
  \draw [->, dotted] (e) edge (d1);
  \draw [->, dotted] (b1) edge (a);
  \draw [->, dotted] (d) edge (e);
\end{tikzpicture}

\medskip

\end{center}
% $$\xymatrix{
%  & x_{n-1} \ar@{.>}[dl] & \cdots \ar@{.>}[l] & x_1
%  \ar@{.>}[l] &
% \\
%  x_n \ar@{-->}[r] & (y-\frac{\mu+1}{2})\xi_n \ar@{->}[r] & p \ar@{->}[r] &
%  (y+\frac{\mu+1}{2})\xi_n \ar@{.>}[r] & x_0 \ar@{.>}[ul]
%  }$$   %AQUI!

Let $\Delta = [[\xi'_0,\xi'_1,\xi'_2]]$. 
Since $\mesh(\Pi) \leq \mu$, there exist $u_i \in \im(\xi'_i)$ such that
$\diam(\{ u_0,u_1,u_2\}) \leq \mu$. Suppose that $u_0 \notin
Z$. Then $u_0 \in [y+\frac{\mu+1}{2},s_n]\xi_n$. Since $u_1 \in
[0,y-\frac{\mu+1}{2}]\xi_n$ and $\xi_n$ is a geodesic, this
contradicts $d(u_0,u_1) \leq \mu$. Hence 
$u_0 \in Z$ and we get
$$d(p,Z) \leq d(p,u_0) \leq d(p,u_2) + d(u_2,u_0) \leq \frac{\mu+1}{2} + \mu =
\frac{3\mu+1}{2}$$
as claimed. Therefore $(X,d)$ is polygon $\frac{\mu+1}{2}$-hyperbolic.
\qed

\section{Virtually free groups}

% Given $\delta \geq 0$, we say that a geodesic metric space $(X,d)$ is:
% \bi
% \item
% {\em polygon $\delta$-hyperbolic} if every geodesic polygon in $(X,d)$ is
% $\delta$-thin; 
% \item
% {\em strongly $\delta$-hyperbolic} if the mesh of every quasi-geodesic
% triangle in $(X,d)$ is $\leq \delta$.
% \ei
% we say that a geodesic metric space is:
% \bi
% \item
% {\em polygon hyperbolic} if it is polygon $\delta$-hyperbolic for
% some $\delta \geq 0$; 
% \item
% {\em strongly hyperbolic} if it is strongly $\delta$-hyperbolic for
% some $\delta \geq 0$.
% \ei

We present next several equivalent geometric
characterizations of virtually free groups. We only have found out that the
equivalence (i) $\iff$ (iv) had been previously proved by Antolin in
\cite{Ant} after we had written our own proof of the theorem, so we
decided to keep our alternative proof of the equivalence
(i) $\iff$ (ii). 

We denote by $(g|h)^A_p$ the Gromov product of $g,h \in G$ with
basepoint $p$ in $(\oo{\Gamma}_A(G),d_A)$.

\bt
\label{fgvf}
Let $G = \langle A \rangle$ be a finitely generated group. Then the following
conditions are equivalent:
\bi
\item[(i)] $G$ is virtually free;
\item[(ii)] $\oo{\Gamma}_A(G)$ is polygon hyperbolic;
% there exists some $\delta \geq 0$ such that
% every geodesic polygon in 
% $\oo{\Gamma}_A(G)$ is $\delta$-thin;
% \beq
% \label{fgvf1}
% \forall x \in [g_0,g_n] \hspace{.7cm} d_A(x,[g_0,g_1] \cup \ldots \cup
% [g_{n-1},g_n]) \leq \delta
% \eeq 
% holds for every geodesic polygon $[[g_0,\ldots,g_n]]$ in 
% $\oo{\Gamma}_A(G)$;
\item[(iii)] there exists some $\delta' \geq 0$ such that
$$(g_0|g_n)^A_p \geq {\rm min}\{ (g_0|g_1)^A_p, \ldots (g_{n-1}|g_n)^A_p
\} -\delta'$$
for all $n \geq 2$ and  $g_0,\ldots, g_n,p \in G$;
\item[(iv)] there exists some $\varepsilon \geq 0$ such that, for all
coterminal geodesic $\xi$ 
and path $\xi'$ in $(\oo{\Gamma}_A(G),d_A)$, 
$${\rm Im}(\xi) \subseteq D_{\varepsilon}({\rm
  Im}(\xi'));$$ 
\item[(v)] ${\rm mesh}(\oo{\Gamma}_A(G),d_A) < +\infty$.
\ei
\et

\proof
(i) $\Rw$ (ii). 
Since $G$ is virtually free and finitely generated, $G$ has 
a finite index free subgroup $F$ of finite rank. By \cite[Proposition
1.11]{GH}, there 
exists a finite generating set $B$ of $F$ such that $\oo{\Gamma}_A(G)$ is
quasi-isometric to $\oo{\Gamma}_B(F)$. 
Let $C$ be a basis of $F$. As mentioned in Section \ref{hg},
$\oo{\Gamma}_B(F)$ is quasi-isometric to $\oo{\Gamma}_C(F)$. Since
$\oo{\Gamma}_C(F)$ 
is a tree, it is trivially polygon $0$-hyperbolic. Therefore, by Lemma
\ref{cuqi}(i), $\oo{\Gamma}_A(G)$ is polygon hyperbolic.

(ii) $\Rw$ (i). Since $\oo{\Gamma}_A(G)$ is polygon hyperbolic, it is
in particular hyperbolic. Thus $G$ is a hyperbolic group. Suppose that
$G$ is not virtually free. By a 
theorem of Bonk and Kleiner \cite[Theorem 1]{BK}, the hyperbolic plane
$\mathbb{H}^2$ admits a quasi-isometric embedding into
$\oo{\Gamma}_A(G)$. In view of Lemma
\ref{cuqi}(i), to reach the required contradiction it suffices to show
that $\mathbb{H}^2$ is not polygon hyperbolic.

We shall use the Poincar\'e disk model of $\mathbb{H}^2$ (see
\cite[p.232]{Gre} for 
details). We can view $\mathbb{H}^2$ as the open disk $\{ z \in
\mathbb{C} \; \big{\vert} \; |z| < 1\}$. Geodesic lines are either
diameters of the 
disk or arcs of a circle whose ends are perpendicular to the disk
boundary. We denote by $d$ the Poincar\'e metric on $\mathbb{H}^2$.

Let $\varepsilon > 0$. Using a diameter and a sufficiently large
number of small circles, 

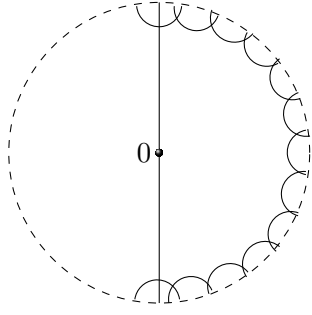
\begin{center}
\begin{tikzpicture}
\draw[dashed] (0,0) circle (2cm);
\draw (0,2) -- (0,-2);
\node at (-0.2,0) {0};
\draw (1.95,0.3) arc (100:260:0.3cm);
\draw (1.85,0.8) arc (110:280:0.3cm);
\draw (1.6,1.25) arc (125:290:0.3cm);
\draw (1.17,1.6) arc (150:300:0.3cm);
\draw (0.7,1.9) arc (160:320:0.32cm);
\draw (0.2,1.95) arc (175:340:0.3cm);
\draw (-0.3,1.95) arc (185:355:0.3cm);
\draw (1.95,-0.25) arc (85:255:0.3cm);
\draw (1.85,-0.8) arc (70:230:0.3cm);
\draw (1.6,-1.25) arc (50:215:0.3cm);
\draw (1.18,-1.59) arc (38:200:0.3cm);
\draw (0.7,-1.83) arc (23:190:0.3cm);
\draw (0.27,-1.95) arc (8:172:0.3cm);
\shade[shading=ball, ball color=black] (0,0) circle (.05);
\end{tikzpicture}
\end{center}

\noindent
we can build a geodesic polygon
$[[\xi_0,\ldots, \xi_n]]$ such that $\im(\xi_n)$ is contained in the
diameter and $|z| > 1-\varepsilon$ for every $z \in \im(\xi_0) \cup
\ldots \cup \im(\xi_{n-1})$. Consider $0 \in \im(\xi_n)$. Then 
$$d(0,z) = \big{\vert}\ln\frac{1-|z|}{1+|z|}\big{\vert}$$
(see \cite[p.248]{Gre}). Since the
euclidean distance from 0 to $\im(\xi_0) \cup
\ldots \cup \im(\xi_{n-1})$ tends to 1 when $\varepsilon$ tends to 0,
it follows that $d(0,\im(\xi_0) \cup
\ldots \cup \im(\xi_{n-1}))$ tends to $+\infty$ when $\varepsilon$
tends to 0. Therefore $\mathbb{H}^2$ is not polygon hyperbolic.

(ii) $\Rw$ (iii) and (ii) $\iff$ (iv) $\iff$ (v). By Theorem \ref{eqg}.

(iii) $\Rw$ (ii). 
If $(X,d)$ is a metric space and $x,y,z,x',y',z' \in X$ are such
$d(x,x'), d(y,y'), d(z,z') \leq \frac{1}{2}$, then
$$\begin{array}{lll}
(x'|y')_{z'}&=&\frac{1}{2}(d(z',x') + d(z',y') - d(x',y'))\\
&\leq&\frac{1}{2}(d(z',z) + d(z,x) + d(x,x') + d(z',z) + d(z,y) +
d(y,y')\\
&&\hspace{1cm}+ d(x,x') - d(x,y) + d(y',y))\\
&\leq&\frac{1}{2}(d(z,x) + d(z,y) - d(x,y) + 3) = (x|y)_z +\frac{3}{2}.
\end{array}$$

Let $n \geq 2$ and  $x_0,\ldots, x_n,q \in
\oo{\Gamma}_AG$. Then there exist $g_0,\ldots, g_n,p \in G$ such that 
$d_A(x_0,g_0), \ldots,$ $d_A(x_n,g_n), d_A(q,p) \leq \frac{1}{2}$. By
the hypothesis and
the inequality proved above, we get
$$\begin{array}{lll}
(x_0|x_n)^A_q&\geq&(g_0|g_n)^A_p-\frac{3}{2} \geq \min\{ (g_0|g_1)^A_p, \ldots
(g_{n-1}|g_n)^A_p\} -\delta' -\frac{3}{2}\\
&\geq&\min\{ (x_0|x_1)^A_q, \ldots (x_{n-1}|x_n)^A_q
\} -\delta' -3.
\end{array}$$
Thus we may apply Theorem \ref{eqg} and obtain condition (ii).
\qed

We remark that the hyperbolicity condition for a geodesic metric space
(all geodesic triangles being $\delta$-thin) is equivalent to all
geodesic $m$-gons being $(m-2)\delta$-thin when we fix $m \geq 4$: for the
direct implication we triangulate the polygon, and the converse
follows from adding trivial geodesics to a triangle. Thus what
distinguishes virtually free groups among arbitrary hyperbolic groups
is the capacity of using the same $\delta$ for all polygons
simultaneously.

With respect to bigons, we must recall a remarkable theorem of
Papasoglu \cite[Theorem 1.4]{Pap}, which proved that a finitely
generated group $G = \langle A \rangle$ is hyperbolic if and only if
there exists some $\delta \geq 0$ such that 
all geodesic bigons in $\oo{\Gamma}_A(G)$ are $\delta$-thin.

As a matter of fact, Papasoglu uses in his statement a synchronous
version of thinness 
(the fellow travel property). We say that two coterminal geodesics
$\xi,\xi':[0,s] \to \oo{\Gamma}_A(G)$ {\em fellow travel} with bound $\delta$
if $d_A(t\xi,t\xi') \leq \delta$ for every $t \in [0,s]$. Clearly,
this implies $\im(\xi) \subseteq D_{\delta}(\im(\xi'))$. 

On the other
hand, suppose that $\im(\xi) \subseteq D_{\varepsilon}(\im(\xi'))$ and
let $t \in [0,s]$. Then $d_A(t\xi,t'\xi') \leq \varepsilon$ for some
$t' \in [0,s]$. It is easy to see that $t' > t + \varepsilon$
contradicts $\xi'$ being a geodesic (we may use $\xi$ and a geodesic
from $t\xi$ to $t'\xi$ to produce a shorter alternative to $\xi'$), and by
symmetry we get $|t'-t| \leq \varepsilon$. Hence
$$d_A(t\xi,t\xi') \leq d_A(t\xi,t'\xi') + d_A(t'\xi',t\xi') \leq 
\varepsilon + |t'-t| \leq 2\varepsilon$$
and so the two geodesics fellow travel with bound $2\varepsilon$. This
confirms that the two notions of thinness are indeed equivalent.

However, we must note that Papasoglu's theorem does not hold for
arbitrary geodesic metric spaces. Indeed, geodesic bigons are
$0$-thin in any space where there exists a unique geodesic connecting
a given pair of points. This is for instance the case of the Euclidean
plane, as noted in \cite{Pap}. Another counterexample would be
$\oo{\Gamma}_A(G)$ when $A = \{ a_1,a_2,\ldots\}$ and $G$ is the group
presented by $\langle A \mid a_n^{2n+1} \; (n \geq 1)\rangle$ (the
free product of cyclic groups of all odd orders).

We end with an application of Theorem \ref{fgvf}.

\bc
\label{aqc}
Let $G$ be a finitely generated virtually free group and let $H \leq
G$. Then the following conditions are equivalent:
\bi
\item[(i)] $H$ is quasi-convex;
\item[(ii)] $H$ is finitely generated.
\ei
\ec

\proof
(i) $\Rw$ (ii). Since every quasi-convex subgroup of a hyperbolic group
is finitely generated \cite[Lemma III.$\Gamma$.3.5]{BH}.

(ii) $\Rw$ (i). Assume that $G = \langle A \rangle$ and $H = \langle B
\rangle$ with $A,B$ finite. Let
$$M = \max\{ d_A(1,b) \mid b \in B\}.$$
By Theorem \ref{fgvf}, there exists some $\varepsilon \geq 0$ such that, for all
coterminal geodesic $\xi$ 
and path $\xi'$ in $(\oo{\Gamma}_A(G),d_A)$, 
$$\im(\xi) \subseteq D_{\varepsilon}(\im(\xi')).$$ 

Let $h,h' \in H$ and let $\eta:[0,s] \to \oo{\Gamma}_A(G)$ be a
geodesic with $0\eta = h$ and $s\eta = h'$. Write $h' = hb_1\ldots
b_n$ for some $b_1,\ldots,b_n \in \wt{B}$. Let  
$\eta'$ be a path in $\oo{\Gamma}_A(G)$ concatenating geodesics with
endpoints $hb_1\ldots b_{i-1}$ and $hb_1\ldots b_i$ for $i =
1,\ldots,n$. Then
$\im(\eta) \subseteq D_{\varepsilon}(\im(\eta'))$.
It is immediate that 
$\im(\eta') \subseteq D_{M}(H)$, hence
$$\im(\eta) \subseteq D_{\varepsilon + M}(H)$$
and $H$ is $(\varepsilon + M)$-quasi-convex.
\qed

\section*{Acknowledgements}

The first author is partially supported by CNPq, PRONEX-Dyn.Syst. and
FAPESB (Brazil).

\smallskip
 
\noindent
The second author acknowledges support from:
\bi
\item
CNPq (Brazil) through a BJT-A grant (process 313768/2013-7);
\item
the European Regional
Development Fund through the programme COMPETE
and the Portuguese Government through FCT (Funda\c c\~ao para a Ci\^encia e a
Tecnologia) under the project\linebreak
PEst-C/MAT/UI0144/2011.
\ei

\end{document}